\documentclass[11pt]{amsart}

\usepackage{amsmath}
\usepackage{amssymb}
\usepackage{mathrsfs}

\newtheorem{theorem}{Theorem}[section]
\newtheorem{claim}[theorem]{Claim}

\newtheorem{lemma}[theorem]{Lemma}

\theoremstyle{definition}
\newtheorem{definition}[theorem]{Definition}

\newtheorem{question}[theorem]{Question}

\theoremstyle{remark}

\newcount\skewfactor
\def\mathunderaccent#1#2 {\let\theaccent#1\skewfactor#2
\mathpalette\putaccentunder}
\def\putaccentunder#1#2{\oalign{$#1#2$\crcr\hidewidth
\vbox to.2ex{\hbox{$#1\skew\skewfactor\theaccent{}$}\vss}\hidewidth}}
\def\name{\mathunderaccent\tilde-3 }

\def\smallbox#1{\leavevmode\thinspace\hbox{\vrule\vtop{\vbox
   {\hrule\kern1pt\hbox{\vphantom{\tt/}\thinspace{\tt#1}\thinspace}}
   \kern1pt\hrule}\vrule}\thinspace}

\newcommand{\cf}{{\rm cf}}
\DeclareMathOperator{\len}{\ell g}
\DeclareMathOperator{\dom}{dom}

\def\qedref#1{$\qed_{\reforiginal{#1}}$}

\title{Dowker filters and Magidor forcing}
\author{Shimon Garti}
\address{Einstein Institute of Mathematics,
 The Hebrew University of Jerusalem,
 Jerusalem 91904, Israel}
\email{shimon.garty@mail.huji.ac.il}

\author{Yair Hayut}
\address{Universit\"at Wien, UZA 1, Institut f\"ur Mathematik, Kurt G\"odel Research Center, Augasse 2-6 1090, Wien, Austria}
\email{yair.hayut@mail.huji.ac.il}

\thanks{The research of the first author was supported by ERC grant 338821.}

\begin{document}
\let\labeloriginal\label
\let\reforiginal\ref

\begin{abstract}
We prove that the existence of a Dowker filter at $\kappa^+$, where $\kappa$ is regular and uncountable, is consistent with $2^\kappa=\kappa^+$.
We also prove the consistency of a Dowker filter at $\mu^+$ where $\mu>\cf(\mu)>\omega$.
This can be forced with $2^\mu=\mu^+$ as well.
\end{abstract}

\maketitle

\newpage

\section{Introduction}

Let $f:\kappa\rightarrow\mathcal{P}(\kappa)$ be a set-mapping.
We shall say that $f$ is \emph{anti-free} iff there are $\alpha<\beta<\kappa$ such that $\alpha\in f(\beta)\wedge\beta\in f(\alpha)$.
Let $\mathscr{F}$ be a filter over $\kappa, e\in{}^\kappa 2$.
A set-mapping $f:\kappa\rightarrow\mathscr{F}$ codes $e$ iff for every $\alpha<\beta<\kappa$, if $\alpha\in f(\beta)\wedge \beta\in f(\alpha)$ then $e(\alpha)=e(\beta)$.
Notice that this is only a partial coding, since if $\alpha\notin f(\beta)$ or $\beta\notin f(\alpha)$ then $f$ gives no information about the relationship between $e(\alpha)$ and $e(\beta)$.

Assume that $\mathscr{F}$ is a filter over $\kappa$.
Is it possible that every set-mapping $f:\kappa\rightarrow\mathscr{F}$ would be anti-free?
A trivial positive answer is given by $\mathscr{F}=\{\kappa\}$, in which case the only function $f:\kappa\rightarrow\mathscr{F}$ is the constant function $f(\alpha)=\kappa$ for every $\alpha\in\kappa$.
Hence for producing an anti-free filter, the simpler the better.

Consider now the coding property.
If $\mathscr{F}=\{\kappa\}$ and $e\in{}^\kappa 2$ is not constant then $e$ cannot be coded by a function $f:\kappa\rightarrow\mathscr{F}$.
Indeed, if $e(\alpha)\neq e(\beta)$ then the only possible $f$ (which is the constant function) fails to code $e$ at $\{\alpha,\beta\}$.
We conclude that $\mathscr{F}$ must be somewhat complicated in order to code $\kappa$-reals, whence we infer that anti-freeness and coding functions are enemies to some extent.

Dowker, \cite{MR0047741}, asked if they can be friends, notwithstanding.
Let $\mathscr{F}$ be a filter over $\kappa$.
Call $\mathscr{F}$ a Dowker filter iff every $f:\kappa\rightarrow\mathscr{F}$ is anti-free and every $e\in{}^\kappa 2$ is coded by some function $f_e:\kappa\rightarrow\mathscr{F}$.
Dowker proved that no Dowker filters live over $\omega_1$ and asked about larger cardinals.
The first positive answer in the environment of the failure of the axiom of choice has been published in \cite{MR973098}.
A few years later, Balogh and Gruenhage forced Dowker filters over $\omega_2$ without violating AC, see \cite{MR1136457}.

An interesting feature of $\omega_2$ is that if $\mathscr{F}$ is a Dowker filter over $\omega_2$ then $\mathscr{F}$ must be uniform.
Namely, if $x\in\mathscr{F}$ then $|x|=\aleph_2$.
Indeed, if $\mathscr{F}$ is a Dowker filter over $\kappa$ then every countable $x\subseteq\kappa$ is $\mathscr{F}$-small.
Likewise, the size of the smallest element of $\mathscr{F}^+$ is strictly less than the minimal cardinality of an element in $\mathscr{F}$.
Hence if $\kappa=\omega_2$ then for some $x\in\mathscr{F}^+$ we have $|x|=\aleph_1$ and every $x\in\mathscr{F}^+$ is uncountable.
It follows that the minimal size of an element of $\mathscr{F}$ is greater than $\aleph_1$ so if $x\in\mathscr{F}$ then $|x|=\aleph_2$.

Non-uniform Dowker filters over every $\kappa>\omega_2$ are obtained immediately from \cite{MR1136457}.
Let $\mathscr{F}$ be a Dowker filter over $\omega_2$, let $\kappa$ be greater than $\aleph_2$ and define $\mathscr{D}=\{A\subseteq\kappa:A\cap\omega_2\in\mathscr{F}\}$.
Such a filter is not uniform, of course.
So the interesting question with respect to cardinals larger than $\omega_2$ takes uniformity as a necessary assumption.
Let us incorporate this property in the formal definition.

\begin{definition}
\label{defdowker} Dowker filters. \newline
Assume that $\kappa\geq\omega_2$ and let $\mathscr{F}$ be a filter over $\kappa$.
We shall say that $\mathscr{F}$ is a Dowker filter iff:
\begin{enumerate}
\item [$(\aleph)$] $\mathscr{F}$ is uniform.
\item [$(\beth)$] $\mathscr{F}$ is anti-free, that is for every $f:\kappa\rightarrow\mathscr{F}$ there are $\alpha<\beta<\kappa$ such that $\alpha\in f(\beta)\wedge\beta\in f(\alpha)$.
\item [$(\gimel)$] $\mathscr{F}$ codes $\kappa$-reals, that is for every $e\in{}^\kappa 2$ there exists $f_e:\kappa\rightarrow\mathscr{F}$ such that $\alpha\in f_e(\beta)\wedge\beta\in f_e(\alpha)\Rightarrow e(\alpha)=e(\beta)$.
\end{enumerate}
\end{definition}

The proof of \cite{MR1136457} has been considerably simplified by Cummings and Morgan in \cite{MR3717964}.
They proved that if $\kappa=\cf(\kappa)>\aleph_0$ and $\mathbb{P}=Add(\kappa,\kappa^{++})$ then there is a Dowker filter over $\kappa^+$ in $V[G]$ whenever $G\subseteq\mathbb{P}$ is $V$-generic.
Notice that $V[G]\models 2^\kappa=\kappa^{++}$ in these models, and thence they raised the following:

\begin{question}
\label{qcm} Is it consistent that $2^\kappa=\kappa^+$ and $\kappa^+$ carries a Dowker filter?
\end{question}

The first possible value of $\kappa=\aleph_0$ which gives a Dowker filter over $\omega_2$ is suggestive, since Dowker filters do not exist over $\omega_1$.
Nevertheless, we shall give a positive answer in the next section by proving that $2^\kappa=\kappa^+$ is consistent with a Dowker filter over $\kappa^+$.
This brings to light another fundamental question, namely whether the existence of a Dowker filter over some successor cardinal is a theorem of ZFC.
An easy negative answer, using the GCH, is excluded by our results.

The previous results concerning Dowker filters were obtained at successors of regular cardinals, and in some sense the forcing constructions are more adapted to this setting.
In particular, Cummings and Morgan raised in \cite{MR3717964} the following:

\begin{question}
\label{qcmsing} Is it consistent that $\mu>\cf(\mu)$ and there is a Dowker filter over $\mu^+$?
\end{question}

If the question is phrased in this way, then the answer is yes.
One can start with a singular cardinal $\mu$ and force with $Add(\theta,\mu^+)$ where $\theta$ is some regular cardinal below $\mu$.
However, it is clear that this is not what the poet meant, and the interesting question is whether one can force a Dowker filter over $\mu^+$ where $\mu$ is a \emph{strong limit} singular cardinal.
Assuming that there is a measurable cardinal $\mu$ with $o(\mu)\geq\omega_1$ in the ground model we will be able to give a positive answer.
Namely, one can force a Dowker filter over $\mu^+$ where $\mu$ is a singular cardinal and $\mu$ is strong limit.

Our notation is standard.
We shall say that $\kappa$ is \emph{strongly regular} iff $\kappa=\kappa^{<\kappa}$.
If $e\in{}^\alpha 2$ then we shall say that $\alpha=\len(e)$.
We use the Jerusalem forcing notation, so $p\leq q$ means that $p$ is weaker than $q$.
The forcing notion $Add(\kappa,\lambda)$ is the usual Cohen forcing for adding $\lambda$-many subsets of $\kappa$ using partial functions from $\lambda\times\kappa$ into $\{0,1\}$ of size less than $\kappa$.
If $\mathscr{F}$ is a filter over $\kappa$ then $\mathcal{I}(\mathscr{F})=\{A\subseteq\kappa:(\kappa-A)\in\mathscr{F}\}$ and $\mathscr{F}^+=\mathcal{P}(\kappa)-\mathcal{I}(\mathscr{F})$.
The elements of $\mathscr{F}^+$ will be called $\mathscr{F}$-positive sets.
If $\mathbb{P}$ is a forcing notion and $\name{\tau}$ is a $\mathbb{P}$-name then a \emph{nice name} for a subset of $\name{\tau}$ is a name of the form $\bigcup\{\{\sigma\}\times A_\sigma:\sigma\in\dom \tau\}$, every $A_\sigma$ being an antichain of $\mathbb{P}$.
Nice names are helpful when one has to reduce the number of names of some object.
The following lemma will be used in the next two sections:

\begin{lemma}
\label{lemcm} Assume that:
\begin{enumerate}
\item [$(a)$] $\mathscr{G}$ is an $\aleph_1$-complete filter over $\mu$.
\item [$(b)$] $C\in\mathscr{G}^+$.
\item [$(c)$] $h:C\rightarrow[{}^{<\mu}{\rm Ord}]^{<\omega}, h(\alpha)\in[{}^{\alpha+1}{\rm Ord}]^{<\omega}$ for every $\alpha\in\mu$.
\end{enumerate}
Then there exist $a\in[{}^\mu{\rm Ord}]^{<\omega}$ and $D\in\mathscr{G}^+, D\subseteq C$ such that for every $f\in{}^\mu{\rm Ord}$ if $f\notin a$ then
\[\{\beta\in D:f\upharpoonright(\beta+1)\in h(\beta) \text { and }\forall g \in a, f \upharpoonright (\beta + 1) \neq g \upharpoonright (\beta + 1)\}\notin\mathscr{G}^+.\]
\end{lemma}

\par\noindent\emph{Proof}. \newline
Assume toward contradiction that the conclusion fails.
For every $n\in\omega$ let $C_n=\{\alpha\in C:|h(\alpha)|=n\}$, so $C=\bigcup_{n\in\omega}C_n$.
By $(a)$ we can choose $n\in\omega$ so that $C_n\in\mathscr{G}^+$.
By finite induction on $i\leq n+1$ we shall define a triple $(a_i,D_i,\eta_i)$ such that:
\begin{enumerate}
\item [$(\aleph)$] $a_i=\{\eta_j:j<i\}\subseteq[{}^\mu{\rm Ord}]^{<\omega}$.
\item [$(\beth)$] $\eta_i\in{}^\mu{\rm Ord}-a_i$.
\item [$(\gimel)$] The set, $D_{i+1}$ which consists of all $\beta \in D_i$ such that $\eta_i\upharpoonright(\beta+1)\in h(\beta)$ and $\forall j < i, \eta_j \upharpoonright (\beta + 1) \neq \eta_i \upharpoonright (\beta + 1)$ is in $\mathscr{G}^{+}$.
\end{enumerate}
Notice that $D_{i+1}$ is determined at the $i$th stage, and we stipulate $D_0=C_n$.
For choosing the elements of the triple at the $(i+1)$th stage we apply the assumption toward contradiction with respect to $(a_i,D_i)$.
This application gives an element $f\in{}^\mu{\rm Ord}, f\notin a_i$ such that $D_{i+1}=\{\beta:\in D_i:f\upharpoonright(\beta+1)\in h(\beta)\}\in\mathscr{G}^+$.
We let $\eta_i=f$ and $a_i=\{\eta_j:j<i\}$.

Choose any ordinal $\beta\in D_{n+1}$ and notice that $\beta\in\bigcap_{i\leq n+1}D_i$.
Hence $\eta_i \upharpoonright (\beta + 1) \in h(\beta)$ for every $i\leq n+1$, so $|h(\beta)|\geq n+1$ by part $(\beth)$ which implies that $i\neq j\Rightarrow \eta_i\upharpoonright (\beta + 1) \neq\eta_j \upharpoonright (\beta + 1)$.
On the other hand $|h(\beta)|=n$ since $\beta\in D_n\subseteq D_0$, a contradiction.

\hfill \qedref{lemcm}

Remark that $h:\mu\rightarrow[{}^{<\mu}{\rm Ord}]^{<\omega}$ but the same statement holds if $h:\mu\rightarrow[Z]^{<\omega}$ for some set $Z$, as in \cite{MR3717964}.
We shall apply this lemma also where $Z=\mu$, in which case the set $a$ will be simply a finite set of ordinals in $\mu$.

Finally, we will use Magidor forcing from \cite{MR0465868} in order to give a positive answer to Question \ref{qcmsing}.
It is not clear whether Prikry forcing can serve in this context, and actually we do not know if a Dowker filter over a successor of a singular cardinal with countable cofinality can be forced.
The essential property of Magidor forcing which we use is a good covering property of new $\omega$-subsets in the Magidor extension.
An explicit proof of the following lemma can be found in \cite[Claim 3.5]{MR3957389}.

\begin{lemma}
\label{lemcovering} Assume that $\kappa$ is measurable, $o(\kappa)\geq\omega_1$ and $\lambda\geq\kappa$. \newline
Let $\mathbb{M}$ be Magidor forcing to make $\kappa>\cf(\kappa)>\omega$, and let $\name{\tau}$ be an $\mathbb{M}$-name of an element in $[\lambda]^{\aleph_0\text{-bd}}$. \newline
Then there are $p\in\mathbb{M}, \theta<\kappa$ and $x\in V$ such that $|x|=\theta$ and $p\Vdash\name{\tau}\subseteq\check{x}$.
\end{lemma}

\hfill \qedref{lemcovering}

The rest of the paper contains two additional sections.
In the first one we shall give a positive answer to Question \ref{qcm} and in the second one we shall try to give a positive answer to Question \ref{qcmsing}.

\newpage 

\section{Dowker filters and the continuum hypothesis}

In this section we prove that the existence of Dowker filters over $\kappa^+$ is consistent with $2^\kappa=\kappa^+$ where $\kappa=\cf(\kappa)>\aleph_0$.
We shall use the ideas of \cite{MR3717964}, which in turn employs the ideas of \cite{MR1136457}.

\begin{theorem}
\label{thmdowkerconthyp} Let $\kappa>\aleph_0$ be strongly regular, let $\mathbb{P}=Add(\kappa,\kappa^{+})$ and let $G\subseteq\mathbb{P}$ be generic over $V$ where $V$ is assumed to satisfy GCH. \newline
Then there is a Dowker filter over $\kappa^+$ in $V[G]$.
\end{theorem}

\par\noindent\emph{Proof}. \newline
Let $\name{e}=(\name{e}_\zeta:\zeta\in\kappa^+)$ be an enumeration of all the nice names of the elements of ${}^{<\kappa^+}2\cap V[G]$ whose length is a successor ordinal.
We will assume that $\ell g(\name{e}_\zeta)$ is decided by the empty condition for every $\zeta\in\kappa^+$.
We require also that for every $\zeta\in\kappa^+$ and every $\delta<\ell g(\name{e}_\zeta)$ there is a unique $\xi\in\kappa^+$ such that $\name{e}_\xi=\name{e}_\zeta\upharpoonright(\delta+1)$.
This is a requirement on the names $\name{e}_\xi,\name{e}_\zeta$, not on the interpretations of them.
It is possible that for some $\xi\neq\xi'$ there will be a condition $p$ such that $p\Vdash\name{e}_\xi=\name{e}_{\xi'}$.

Rather than $Add(\kappa,\kappa^+)$ we shall use the forcing notion $\mathbb{P}$ defined as follows.
A condition $p\in\mathbb{P}$ is a partial function from $\kappa^+\times\kappa^+$ into $\{0,1\}$ so that $|p|<\kappa$.
If $p,q\in\mathbb{P}$ then $p\leq_{\mathbb{P}}q$ iff $p\subseteq q$.
One verifies easily that $\mathbb{P}$ and $Add(\kappa,\kappa^+)$ are isomorphic.

If $G\subseteq\mathbb{P}$ is $V$-generic then $g=\bigcup G$ is a function from $\kappa^+\times\kappa^+$ into $\{0,1\}$.
For each $\xi\in\kappa^+$ let $g_\xi = g\upharpoonright\{\xi\}\times\kappa^+$, so $g_\xi:\kappa^+\rightarrow 2$.
If $y\subseteq\kappa^+\times\kappa^+$ then $\pi_0(y)=\{\alpha\in\kappa^+:\exists\beta,(\alpha,\beta)\in y\}$ and $\pi_1(y)=\{\beta\in\kappa^+:\exists\alpha,(\alpha,\beta)\in y\}$.
We shall use these conventions mainly with respect to the domain of conditions in $\mathbb{P}$.
For every $\zeta\in\kappa^+$ we define a set $A_\zeta\subseteq\kappa^+$ by describing its characteristic function $\chi_{A_\zeta}$ as follows.
If $\beta<\len(e_\zeta)=\eta_\zeta+1$ and $e_\zeta\upharpoonright(\beta+1)=e_\xi$ and $g_\zeta(\beta)=g_\xi(\eta_\zeta)=1$ and $e_\zeta(\beta)\neq e_\zeta(\eta_\zeta)$ then let $\chi_{A_\zeta}(\beta)=0$.
In all other cases let $\chi_{A_\zeta}(\beta)=g_\zeta(\beta)$.

The idea is much simpler than it looks like.
Our default for $\chi_{A_\zeta}$ is the function $g_\zeta(\beta)$, but we rewrite the values of $g_\zeta(\beta)$ in case where we anticipate troubles with the coding property.
This happens, basically, if $g_\zeta(\beta)=g_\zeta(\eta_\zeta)=1$ and concomitantly $e_\zeta(\beta)\neq e_\zeta(\eta_\zeta)$.
Now by removing $\beta$ from $A_\zeta$, to wit by setting $\chi_{A_\zeta}(\beta)=0$, we need not worry about coding at $e_\zeta(\beta)$.
Notice that the length of each $e_\zeta$ and the identity of the ordinal $\xi$ so that $e_\xi=e_\zeta\upharpoonright(\beta+1)$ are computed on names and this can be rendered in $V$, but the values of $g_\zeta(\beta),g_\xi(\eta_\zeta)$ and the inequality $e_\zeta(\beta)\neq e_\zeta(\eta_\zeta)$ are done in $V[G]$.

The sets of the form $A_\zeta$ will generate our Dowker filter.
We shall also define another filter $\mathscr{G}$ from these sets, needed for the proof.
Let $\mathscr{D}$ be the filter obtained by taking all the finite intersections of $A_\zeta$s and closing upwards.
Let $\mathscr{G}$ be the filter obtained by $(<\kappa)$-intersections of $A_\zeta$s, clubs of $\kappa^+$ and the set $S^{\kappa^+}_\kappa$.
It follows that $\mathscr{G}$ is a $\kappa$-complete filter over $\kappa^+$ and it extends the club filter of $\kappa^+$ restricted to $S^{\kappa^+}_\kappa$.
The fact that $\mathscr{D}$ and $\mathscr{G}$ are proper filters (that is, $\varnothing\notin\mathscr{D}$ and $\varnothing\notin\mathscr{G}$) follows essentially from the fact that $|p|<\kappa$ for every $p\in\mathbb{P}$.
A detailed argument appears in \cite[Lemma 4.7]{MR3717964} and \cite[Corollary 4.9]{MR3717964}.

Observe that $\mathscr{D}$ will satisfy Definition \ref{defdowker}($\gimel$) by the definition of the sets $A_\zeta$.
Likewise, each $A_\zeta$ is of size $\kappa^+$, hence $\mathscr{D}$ is uniform.
To see that $\mathscr{D}$ satisfies \ref{defdowker}($\gimel$), assume that $\name{d}$ is a name of a function from $\kappa^+$ into $\{0,1\}$, so it is interpreted as an element of ${}^{\kappa^+}2\cap V[G]$.
Without loss of generality, $\name{d}$ is a nice name.
By our assumption on the enumeration $\name{e}$, for every $\alpha\in\kappa^+$ there is a unique $\zeta(\alpha)\in\kappa^+$ such that $\name{d}\upharpoonright(\alpha+1)=\name{e}_{\zeta(\alpha)}$.
Notice that $\zeta(\alpha)$ can be computed from the name $\name{d}$ in the ground model and define $F(\alpha)=A_{\zeta(\alpha)}$ for every $\alpha\in\kappa^+$.

Assume now that $\alpha<\beta<\kappa^+$ and $\alpha\in F(\beta)\wedge\beta\in F(\alpha)$.
It follows that $g_{\zeta(\alpha)}(\beta)=g_{\zeta(\beta)}(\alpha)=1$.
Recall that $\alpha<\beta+1=\ell g(e_{\zeta(\beta)})$ and by our choice $e_{\zeta(\beta)}\upharpoonright(\alpha+1)=e_{\zeta(\alpha)}$.
Since $g_{\zeta(\alpha)}(\beta)=g_{\zeta(\beta)}(\alpha)=1$ we conclude that $e_{\zeta(\beta)}(\beta)=e_{\zeta(\beta)}(\alpha)$, so $\name{d}$ is coded correctly in $V[G]$ as required.

The burden of the proof is to show that \ref{defdowker}($\beth$) holds in the generic extension.
Since the sets $A_\zeta$ form a base for $\mathscr{D}$, suffice it to deal with functions from $\kappa^+$ into intersections of $A_\zeta$s.
These can be represented by functions from $\kappa^+$ into $[\kappa^+]^{<\omega}$, where for every $\alpha\in\kappa^+$ the function assigns a finite set $\{\zeta_1,\ldots,\zeta_n\}$ and gives rise to $\bigcap_{1\leq i\leq n}A_{\zeta_i}$.
Proving anti-freeness for the intersections $\bigcap_{1\leq i\leq n}A_{\zeta_i}$ will accomplish the proof of Dowkerness.

Suppose that $\name{t}:\kappa^+\rightarrow[\kappa^+]^{<\omega}$.
Define $h:\kappa^+\rightarrow[{}^{<\kappa^+}2]^{<\omega}$ by letting $h(\alpha)=\{e_{\xi_1},\ldots,e_{\xi_n}\}$ whenever $\name{t}(\alpha)=\{\xi_1,\ldots,\xi_n\}$.
We apply Lemma \ref{lemcm} to the filter $\mathscr{G}$ and the function $h$ (where $\mu=\kappa^+$ and $C=\kappa^+$).
By the conclusion of the lemma we infer that there are a set $D\in\mathscr{G}^+$ and a finite set $a\subseteq{}^{\kappa^+}2$ such that:
\begin{enumerate}
\item [$(i)$] For some fixed $n\in\omega$ and for every $\beta\in D$ we have $|h(\beta)|=n$.
\item [$(ii)$] For every $f:\kappa^+\rightarrow 2$ if $f\notin a$ then $\{\beta\in D:f\upharpoonright(\beta+1)\in h(\beta)\}\notin\mathscr{G}^+$.
\end{enumerate}
Note that we can remove the assumption that $f \upharpoonright (\beta + 1) \neq g \upharpoonright (\beta + 1)$ for all $g \in a$ that appears in the formulation of the lemma, since the filter $\mathscr{G}$ contains all the tails of $\kappa^+$, and thus, by removing a bounded initial segment from $D$, we may assume that the initial segments of elements of $a$ are different from the initial segments of $f$.

Let $a=\{f_0,\ldots,f_{k-1}\}$.
There are only finitely many possibilities for the sequence $(f_0(\alpha),\ldots,f_{k-1}(\alpha))$.
Hence by shrinking $D$ further we may assume that the sequence $(f_0(\alpha),\ldots,f_{k-1}(\alpha))$ is the same fixed sequence for every $\alpha\in D$.
We apply the lemma once again to obtain a finite set of ordinals $a'\subseteq\kappa^+$ such that if $\gamma\notin a'$ then $\{\beta\in D:\gamma\in\name{t}(\beta)\}\notin\mathscr{G}^+$.
For this, see the paragraph immediately after Lemma \ref{lemcm}.

Let us work in $V$. Let $\chi$ be a sufficiently large regular cardinal.
Choose $M\prec\mathcal{H}(\chi)$ which enjoys the following properties:
\begin{enumerate}
\item [$(\aleph)$] $|M|=\kappa$.
\item [$(\beth)$] ${}^{<\kappa}M\subseteq M$.
\item [$(\gimel)$] $a,a'\in M$.
\item [$(\daleth)$] $\beta=M\cap\kappa^+\in\kappa^+$.
\end{enumerate}
Let $p_0\in\mathbb{P}$ be a condition which forces a value to $n$ and determines the elements of $a'$.
Choose a condition $q\in\mathbb{P}$ such that $p_0\leq q$ and $q\Vdash\check{\beta}\in\name{D}$.
To obtain this model $M$ and the ordinal $\beta$, let $(M_i:i\in\kappa^+)$ be an increasing continuous sequence of elementary submodels of $\mathcal{H}(\chi)$ satisfying the above requirements and notice that $E=\{M_i\cap\kappa^+:i\in\kappa^+\}$ is a club of $\kappa^+$.
Recall that $\mathscr{G}$ extends the club filter of $\kappa^+$ and $D\in\mathscr{G}^+$, hence for some $i\in\kappa^+$ we will have a condition $q\geq p_0$ such that $q\Vdash M_i\cap\kappa^+\in E\cap D$.
Set $M=M_i$ and $\beta=M_i\cap\kappa^+$.
Observe that by the same reasoning, we can pick $M_i$ in a way that $q\Vdash\beta\in\bigcap_{\zeta\in a'}A_\zeta$ since $\bigcap_{\zeta\in a'}A_\zeta\in\mathscr{G}$.

By extending $q$ if needed we can assume that $q$ forces a value to $\name{t}(\beta)$, that is $q\Vdash\name{t}(\beta)=\{\xi_1,\ldots,\xi_n\}$.
Likewise, $q\Vdash\len(e_{\xi_i})=\gamma_i$ for every $1\leq i\leq n$.
Observe that $\cf(\beta)=\kappa$ since $\beta=M\cap\kappa^+$ and ${}^{<\kappa}M\subseteq M$.
Hence requirement $(\beth)$ implies that $q\cap M\in M$ bearing in mind that $|q|<\kappa$.

Working within $M$ we are asking whether there is a condition $r\in M$ such that $q\cap M\leq r$ and for some ordinal $\alpha\in\kappa^+$ we can extend $r$ to a condition which forces $\alpha\in A_{t(\beta)}\wedge \beta\in A_{t(\alpha)}$.
More specifically, we wish to find an ordinal $\alpha>\sup(q\cap M), \alpha\in D\cap(\bigcap_{\xi_i\in M}A_{\xi_i}-\bigcup_{\gamma_i\in M}\gamma_i)$ and a condition $r\geq q\cap M$ such that $r\Vdash\alpha\in A_{t(\beta)}\wedge\beta\in A_{t(\alpha)}$.

Let us assume that there are no such $r$ and $\alpha$. Let $D'$ be the set of all $\alpha \in D\cap(\bigcap_{\xi_i\in M}A_{\xi_i}-\bigcup_{\gamma_i\in M}\gamma_i)$ such that $\name{t}(\alpha)-a'$ is disjoint from $\pi_0(\dom q) \cap M$.
Since $|q| < \kappa$, we remove $<\kappa$ many null sets from a positive set and thus $D'$ is forced to be positive and in particular, unbounded in $\kappa^{+}$. Since $D' \in M$, we can find many ordinals $\alpha$ and conditions $r \geq q \cap M$ such that $\alpha, r\in M$ and $r \Vdash \alpha \in D'$. By taking a stronger condition if needed we may assume that $r$ decides the value of $t(\alpha) = \{\zeta_1, \dots, \zeta_n\}$ and the lengths of $e_{\zeta_i}$ for $1 \leq i \leq n$. Note that all those ordinals are going to be in $M$ and in particular, below $\beta$.

By the definition of $A_{\zeta_i}$, $\beta \in A_{\zeta_i}$ if and only if $g_{\zeta_i}(\beta) = 1$. Let us look at the condition $q \cup r$. In order for the pair $(\zeta_i, \beta)$ to be in the domain of this condition, it has to be in $q$.
Indeed, $(\zeta_i,\beta)\notin\dom  r$ since $r\in M$ and $\beta\notin M$ so necessarily $(\zeta_i,\beta)\in\dom q$.
Note that $\zeta_i \in a'$. Otherwise, it contradicts our choice of $D'$. In this case, $g_{\zeta_i}(\beta) = 1$, by our choice of $\beta$. In case that $(\zeta_i,\beta)\notin \dom q$, we can extend $q \cup r$ to a condition $q_1$ by adding $(\zeta_i, \beta)$ to the domain with value $1$. So in any case we found a condition, $q_1$, which forces that $\beta$ is in $A_{t(\alpha)} = \bigcap_i A_{\zeta_i}$.

Let us try to obtain also the other requirement, $\alpha \in A_{t(\beta)}$. Let $\xi_i \in t(\beta)$, and $\gamma_i = \len(e_{\xi_i})$. If $\gamma_i < \beta$, it is in $M$ and thus $\alpha > \gamma_i$. In this case, in order to get $\alpha$ to be in $A_{\xi_i}$, we need to have $g_{\xi_i}(\alpha) = 1$. If further $\xi_i \in M$ then this is part of the definition of $D'$.
Otherwise, $\xi_i \notin M$ and in particular $(\xi_i, \alpha) \notin \dom q \cup \dom r$ since $\alpha \notin \pi_1(\dom  q)$ and $\xi_i \notin \pi_0(\dom  r)$.
Moving to $q_1$ did not change this fact, and thus we can extend $q_1$ to $q_2$ by adding a single pair into the domain of the condition, sending $(\xi_i, \alpha)$ to $1$.

Thus, we are left with the case $\gamma_i \geq \beta$. This implies that $\xi_i$ is not in $M$ as otherwise, by elementarity, the length of $e_{\xi_i}$ would be inside $M$. Thus, the pair $(\xi_i, \alpha)$ is not in the domain of $q_2$ so let $q_3$ be a condition that sends this pair to $1$. Recall that in the definition of $A_{\xi_i}$, in order for $q_3$ not to force $\alpha$ to be in $A_{\xi_i}$ the following must hold: let $e_{\rho} = e_{\xi_i} \upharpoonright(\alpha + 1)$, then $g_\rho(\gamma_i) = g_{\xi_i}(\alpha) = 1$ and $e_\rho(\alpha) = e_{\xi_1}(\alpha) \neq e_{\xi_i}(\gamma_i)$. Since $\gamma_i \notin M$, in order for $q_3$ to force $g_{\rho}(\gamma_i) = 1$, either $\rho \in \pi_0(\dom q)$ or $\gamma_i = \beta$ and $\rho = \zeta_j$ for some $j$. The first case can be avoided easily, since $\len(e_\rho) = \alpha + 1$ and thus there are at most $<\kappa$ many possible values of $\alpha$ for which this can occur.

Note that if $e_{\xi_i}$ is a restriction of $f_k$ for some $f_k\in a$, then since we stabilized the values of $f_k$ on $D'$ and they equal to the values which are obtained in $\beta$, $e_{\xi_1}(\alpha) \neq e_{\xi_i}(\gamma_i)$ is impossible.
Thus, $e_{\xi_i}$ is not a restriction of an element in $a$.
We conclude that this is the only possible obstacle. Let us summarize the discussion so far:

\begin{claim}
\label{clmobstacle}
In any generic extension, for all sufficiently large $\alpha \in D' \cap M$, there are $i, j \leq n$ such that $e_{\xi_i} \upharpoonright (\alpha + 1) = e_{\zeta_j}$ and $e_{\xi_i}$ is not a restriction of an element in $a$.
\end{claim}

Let us work now in the generic extension $V[G]$. By standard arguments, $M[G \cap M] \prec \mathcal{H}(\chi)[G]$. In $V[G]$, one can partition an end segment of $D' \cap M$ into $\leq n^2$ many pieces according to the values of $i$ and $j$ in Claim \ref{clmobstacle}. In particular, in $M[G \cap M]$, for every finite collection of ordinals in $D' \cap M$, $F$, there is a partition into $n^2$ many parts $F = \bigcup F_{i,j}$, such that for every $\alpha, \beta \in F_{i,j}$, if $\alpha < \beta$ then $e_{\zeta_j^\alpha}$ is an initial segment of length $\alpha + 1$ of $e_{\zeta_j^\beta}$.

By elementarity, the same holds in $V[G]$. Thus, in $V[G]$, every finite subset of $D'$ can be colored by $n^2$ colors such that if $\alpha < \beta$ are both colored in the color $(i,j)$ then $e_{\zeta_j}^\alpha$ is an initial segment of $e_{\zeta_j}^\beta$ (with the right lengths), and $e_{\zeta_j^\alpha}$ is not an initial segment of any $f\in a$.
By compactness, we conclude that there is a coloring of all $D'$ into $n^2$ many colors with this property.
Pick a positive homogeneous set, say $H_{i,j}$ and take $f = \bigcup_{\alpha \in H_{i,j}} e_{\zeta_j^\alpha}$. Then $f\notin a$, which is a contradiction to Lemma \ref{lemcm}.
We conclude, therefore, that for some $r\in M, \alpha\in {\rm Ord}$ it is true that $r\Vdash\alpha\in A_{\name{t}(\beta)} \wedge \beta\in A_{\name{t}(\alpha)}$.
Hence $\mathscr{D}$ is a Dowker filter and we are done.

\hfill \qedref{thmdowkerconthyp}

We conclude this section by mentioning the original question of Dowker, which remains open:

\begin{question}
\label{qdowkerzfc} Does there exist a successor cardinal $\kappa^+$ which carries a Dowker filter in ZFC?
\end{question}

\newpage 

\section{Successor of a singular cardinal}

The theorem in the previous section gives a positive answer to Question \ref{qcm} when $\kappa=\cf(\kappa)>\aleph_0, 2^\kappa=\kappa^+$ and a Dowker filter is forced over $\kappa^+$.
The reason for $\kappa>\aleph_0$ is that we need the filter $\mathscr{G}$ within the proof to be $\omega_1$-complete.
This plays a key-role in the proof, and if $\kappa=\aleph_0$ or more generally if $\cf(\kappa)=\omega$ then this essential point becomes problematic.
The concrete case of $\aleph_1$ shows that this issue is not merely a technical problem.
For this reason, a simple attempt to force a Dowker filter over $\kappa^+$ when $\kappa$ is measurable and then to add a Prikry sequence to $\kappa$ will probably fail.
However, this obstacle is also suggestive, since we can singularize a large cardinal with Magidor forcing.

\begin{theorem}
\label{thmmagforcing} Suppose that $\mu$ is a measurable cardinal and $o(\mu)\geq\omega_1$. \newline
Then one can force $\mu>\cf(\mu), \mu$ is strong limit along with the existence of a Dowker filter over $\mu^+$.
\end{theorem}

\par\noindent\emph{Proof}. \newline
Let $\mathbb{P}$ be an iteration with Easton support of $Add(\delta,\delta^+)$ at every inaccessible cardinal $\delta\leq\mu$, and let $G\subseteq\mathbb{P}$ be generic over $V$.
Remark that $\mu$ is still measurable and $o(\mu)\geq\omega_1$ in $V[G]$.
Working in $V[G]$, let $\mathbb{M}$ be Magidor forcing to make $\mu>\cf(\mu)=\omega_1$ and let $H\subseteq\mathbb{M}$ be generic over $V[G]$.
We claim that in $V[G\ast H]$ there is a Dowker filter over $\mu^+$.

Let $\name{e}=(\name{e}_\zeta:\zeta\in\mu^+)$ be an enumeration of all the nice names of the elements of ${}^{<\mu^+}2\cap V[G\ast H]$ whose length is a successor ordinal.
We may assume that $\ell g(\name{e}_\zeta)$ is decided by the empty condition for every $\zeta\in\mu^+$.
As in the previous section we let $g=\bigcup G$ and $g_\zeta=g\upharpoonright\{\zeta\}\times\mu^+$ for every $\zeta\in\mu^+$, so each $g_\zeta$ is a function from $\mu^+$ to $2$.

We define $A_\zeta\subseteq\mu^+$ as done in Theorem \ref{thmdowkerconthyp}, for every $\zeta\in\mu^+$.
Explicitly, if $\beta<\ell g(e_\zeta)=\eta_\zeta+1, e_\xi=e_\zeta\upharpoonright(\beta+1), g_\zeta(\beta)=g_\xi(\eta_\zeta)=1$ and $e_\zeta(\beta)\neq e_\zeta(\eta_\zeta)$ then $\chi_{A_\zeta}(\beta)=0$.
Otherwise, $\chi_{A_\zeta}(\beta)=g_\zeta(\beta)$ and of course $A_\zeta$ is the set whose $\chi_{A_\zeta}$ is the characteristic function.

Working in $V[G]$, let $\mathscr{G}'$ be the filter generated by the sets $A_\zeta$ for $\zeta\in\mu^+$, the club filter of $\mu^+$ and $S^{\mu^+}_\mu$.
Let $\mathscr{G}$ be the filter generated by $\mathscr{G}'$ in $V[G\ast H]$.
We shall prove in Lemma \ref{lemcompleteness} below that $\mathscr{G}$ is $\aleph_1$-complete in $V[G\ast H]$ and notice that it contains the $\aleph_1$-complete filter generated by the $A_\zeta$s.
Therefore, we will be able to apply Lemma \ref{lemcm} to $\mathscr{G}$.

Let $\mathscr{D}$ be the filter generated by finite intersections of $A_\zeta$s.
Notice that $\mathscr{D}$ is uniform and satisfies \ref{defdowker}($\gimel$) by the properties of each $A_\zeta$.
The argument is exactly as in the previous section.
It remains to prove that $\mathscr{D}$ satisfies \ref{defdowker}($\beth$) as well.
Assume that $\name{t}:\mu^+\rightarrow[\mu^+]^{<\omega}$ so $\name{t}(\alpha)\in[\mu^+]^{<\omega}$ for each $\alpha\in\mu^+$, and let $A_{\name{t}(\alpha)}=\bigcap\{A_\zeta:\zeta\in\name{t}(\alpha)\}$.
We wish to prove that for every $\name{t}$ there are $\alpha<\beta<\mu^+$ and a condition $r$ so that $r\Vdash\alpha\in A_{\name{t}(\beta)} \wedge \beta\in A_{\name{t}(\alpha)}$.

Given $\name{t}$ we define $h=h(\name{t}):\mu^+\rightarrow[{}^{<\mu^+}2]^{<\omega}$ by $h(\alpha)=\{e_{\xi_1},\ldots,e_{\xi_n}\}$ where $\name{t}(\alpha)=\{\xi_1,\ldots,\xi_n\}$ and let $C=\mu^+$.
Applying Lemma \ref{lemcm} to the triple $(\mathscr{G},h,C)$ we obtain $D\subseteq\mu^+,D\in\mathscr{G}^+$ and a finite set $a$ of functions as guaranteed by the lemma.
For every $\beta\in D$ and every $q\in\mathbb{P}$ we define:
$$
R_{q\beta} = \{\gamma\in\beta:(\gamma,\beta)\in\dom q, q(\gamma,\beta)=0\}.
$$
Working in $V[G\upharpoonright\mu]$ we fix a sufficiently large regular cardinal $\chi$ and an increasing continuous sequence $(M_i:i\in\mu^+)$ of elementary submodels of $\mathcal{H}(\chi)$ such that $|M_i|=\mu, \mu+1\subseteq M_i$ for every $i\in\mu^+$ and each $M_{i+1}$ is $<\mu$-closed.
Notice that if $\cf(i)=\mu$ then $M_i$ will be $<\mu$-closed as well.
We shall prove in Lemma \ref{lemrqbeta} that there exist an ordinal $\beta\in\mu^+$ and a condition $(q,p)\in\mathbb{P}\ast\mathbb{M}$ with the following properties:
\begin{enumerate}
\item [$(a)$] $\beta=M_i\cap\mu^+$ for some $i\in\mu^+$ and $\cf(\beta)=\mu$.
\item [$(b)$] $(q,p)\Vdash\check{\beta}\in\name{D}$.
\item [$(c)$] $(q,p)\Vdash\name{t}(\alpha)\cap R_{q\beta}=\varnothing$ for a $\mathscr{G}$-positive set of $\alpha\in D$.
\end{enumerate}
By extending $(q,p)$ if needed we may assume that $(q,p)$ forces a value to $\name{t}(\beta)$, determines the size of $a$ and forces the value of $\ell g(e_\zeta)$ for every $\zeta\in\name{t}(\beta)$.
Namely, $(q,p)\Vdash\name{t}(\beta)=\{\xi_1,\ldots,\xi_n\}$ and $\ell g(e_{\xi_i})=\gamma_i$ for every $1\leq i\leq n$.

Let $M=M_i$ where $i\in\mu^+$ is the ordinal provided by part $(a)$.
Notice that $\cf(i)=\mu$ and hence ${}^{<\mu}M\subseteq M$.
Every $q\in\mathbb{P}$ satisfies $|\dom q|<\mu$ hence $|R_{q\beta}|<\mu$ as well, and therefore $R_{q\beta}\in M$.
Likewise, the name of the set of all $\alpha\in D$ such that $\name{t}(\alpha)\cap R_{q\beta}=\varnothing$ belongs to $M$ since it is definable in $M$.
For every $1\leq i\leq n$ let $E_i$ be $\{\alpha<\beta:\exists\rho, \name{e}_\rho=\name{e}_{\xi_i}\upharpoonright(\alpha+1), (\rho,\gamma_i)\in\dom q\}$.
Again, $|E_i|<\mu$ and hence the set $E=\bigcup\{E_i:1\leq i\leq n\}$ is of size less than $\mu$, and consequently $E\in M$.

Define $D'_0=\bigcap\{A_{\xi_j}:\xi_j\in M\}-(\bigcup_{\gamma_j\in M}\gamma_j\cup E)$, so $D'_0\in M$.
Finally, let $D'$ be $D\cap D'_0$ intersected with the set of $\alpha\in D$ for which $\name{t}(\alpha)\cap R_{q\beta}=\varnothing$.
Notice that $D'\in M$ and it is forced (by the rest of $\mathbb{P}\ast\mathbb{M}$) to be in $\mathscr{G}^+$.
We work now in $V[G\ast H]$ and we assume that $(q,p)\in G\ast H$.
Observe that $M[(G\ast H)\cap M]\prec\mathcal{H}(\chi)[G\ast H]$, since $\mu\subseteq M$ and $\mathbb{P}\ast\mathbb{M}$ is $\mu^+$-cc.
As noted above, in $M[G\ast H]$ we know that the set of $\alpha\in D$ such that $\name{t}(\alpha)\cap R_{q\beta}=\varnothing$ is in $\mathscr{G}^+$, so in particular it is unbounded in $\mu^+$.

For every $\alpha$ in this positive set we choose a condition $(r_\alpha,p_\alpha)\in G\ast H$ so that $(r_\alpha,p_\alpha)$ forces $\alpha$ to be in this set and it forces the values of $\name{t}(\alpha)$ and the lengths of its elements.
Applying the pigeonhole principle, there are a condition $r\in\mathbb{P}$ and a fixed stem $s_\star$ such that for some $B_{s_\star}\subseteq\mu^+$, $r$ forces that $|B_{s_\star}|=\mu^+$ and if $\alpha\in B_{s_\star}$ then $r_\alpha \geq r$ and the stem of $p_\alpha$ is $s_\star$.
By the directness of $H$ we may assume that the stem of $p$ is $s_\star$ as well.

We claim that there are a condition $(q',p')\in\mathbb{P}\ast\mathbb{M}$ and an ordinal $\alpha\in B_{s_\star}\cap\beta$ such that $(q,p)\leq_{\mathbb{P}\ast\mathbb{M}}(q',p')$ and $(q',p')\Vdash\alpha\in A_{\name{t}(\beta)} \wedge \beta\in A_{\name{t}(\alpha)}$.
The easy part is to force $\beta\in A_{\name{t}(\alpha)}$.
To this end, we choose a condition $(r,p^+)\in M$ and some $\alpha\in M$ such that $\alpha$ is forced to be in $B_{s_\star}, s_\star$ is the stem of $p^+$ and $q\cap M\leq r$.
We indicate that the Cohen part of the condition forces $p\parallel p^+$ since $s_\star$ is the Magidor stem of both conditions.
Hence $r\cup q$ also forces this fact.

Now $(r,p^+)\in M$ and hence $\name{t}(\alpha)\in M$ and for every $\xi\in\name{t}(\alpha)$ we have $\ell g(e_\xi)\in M$ so $\ell g(e_\xi)<\beta$.
By a finite sequence of extensions of $r\cup q$ we can make sure that $\beta$ is forced into $A_\xi$ for every $\xi\in\name{t}(\alpha)$.
To see this, observe first that $\beta$ is absent from the ordinals mentioned in $r$ since $r\in M$ and $\beta\notin M$.
For the $q$ part either $(\xi,\beta)\in\dom q$, in which case $q\Vdash\beta\in A_\xi$ since $\alpha\in D'$, or $(\xi,\beta)\notin\dom q$ and we can add the triple $(\xi,\beta,1)$ to $q$.
As $\name{t}(\alpha)$ is finite we can extend $r\cup q$ and obtain $q'$ and $p'\geq p,p^+$ so that $(q',p')\Vdash\beta\in A_{\name{t}(\alpha)}$.

The harder part is to ensure that $(q',p')\Vdash\alpha\in A_{\name{t}(\beta)}$.
Now if $\xi\in\name{t}(\beta)\cap M$ then $\alpha\in A_\xi$ since $\alpha\in A\subseteq D'\subseteq\bigcap\{A_\zeta:\zeta\in M\}$.
Assume, therefore, that $\xi\in\name{t}(\beta)$ and $\xi\notin M$.
Let $\gamma$ be such that $\ell g(e_\xi)=\gamma+1$.
We distinguish three cases:

\par\noindent\emph{Case 1}: $\gamma<\beta$.

In this case necessarily $\alpha>\gamma$ since $\alpha\in D'$ and by the definition of $D'$.
We claim that one can add the triple $(\xi,\alpha,1)$ to the finite sequence of extenstions of $r\cup q$ rendered to obtain $\beta\in A_{\name{t}(\alpha)}$, which is exactly what we need.
For this, observe that $\xi\notin\pi_0(\dom r)$ since $r\in M$ and $\xi\notin M$.
Likewise, $\alpha\notin\pi_1(\dom q)$ since $\alpha\in D'$ and hence $\alpha\notin E$.
Finally, the finite sequence of extensions of $r\cup q$ needed to force $\beta\in A_{\name{t}(\alpha)}$ involve $\beta$ at each step, and here $\alpha\neq\beta$ so there is no problem with our triple.

\par\noindent\emph{Case 2}: $\gamma>\beta$.

Let $\rho$ be such that $e_\rho=e_\xi\upharpoonright(\alpha+1)$.
Since $\alpha\notin E$ we see that $(\rho,\gamma)\notin\dom q$.
Since $\gamma>\beta$ and $r\in M$ we see that $(\rho,\gamma)\notin\dom r$.
Hence we can add the triple $(\rho,\gamma,0)$ to the extension of $r\cup q$ obtained so far.
By doing so we conclude that our extension of $r\cup q$ does not force $\alpha\notin A_\xi$, hence we can add $(\xi,\alpha,1)$ to this extension.

\par\noindent\emph{Case 3}: $\gamma=\beta$.

In this case, one of the extensions of $r\cup q$ (made for ensuring that $\beta$ is forced into $A_{\name{t}(\alpha)}$) might be problematic when trying to force $\alpha$ into $A_{\name{t}(\beta)}$.
Let us analyze the situation and see how the problem can be avoided.
We shall need the concept of a \emph{thread}.
A set $F=\{f_0, \dots, f_{k-1}\}$ of functions from $\mu^+$ into $2$ is a thread covering for $A$, where $A$ is an unbounded subset of $\mu^+$, iff for every $\alpha\in A$ there is $\zeta\in\name{t}(\alpha)$ and $i < k$ such that $e_\zeta=f_i\upharpoonright(\alpha+1)$. We will always assume that our threads are \emph{minimal}, namely that for every $i < k$, $\{f_j \mid j \neq i\}$ is not a thread.

It may happen that for some $\alpha\in B_{s_\star}$ we will have $\zeta\in\name{t}(\alpha),\xi\in\name{t}(\beta), e_\zeta=e_\xi\upharpoonright(\alpha+1)$ and $\ell g(e_\xi)=\beta+1$.
In such cases $e_\zeta\lhd e_\xi$ (recall that $\alpha<\beta$) and maybe we cannot add $(\xi,\alpha,1)$ to our condition.
We can try a different $\alpha\in B_{s_\star}$, but perhaps the same problem occurs at every $\alpha\in B_{s_\star}$.
In such cases there is a thread covering $F$ for $B_{s_\star}$, using the compactness argument from the proof of Theorem \ref{thmdowkerconthyp}. More precisely, $|F| \leq n$.
Lest $F\cap a \neq \emptyset$ there will be no problem, by the definition of $D'$.
Therefore, our purpose is to show that there exists a stem $s$ such that if there is a thread covering $F = \{f_0, \dots, f_{k-1}\}$, $k \leq n$ for $B_s$ then necessarily $f_i\in a$ for some $i$.

Assume towards contradiction that for every Magidor-stem $s$ there is a thread covering $F_s$ for the set $B_s$ and $F_s \cap a = \emptyset$, $|F_s| \leq n$.
Notice that if $s$ is fixed then there are at most finitely many minimal thread coverings for $B_s$ of size $\leq n$. Recall that $n$ is the size of $\name{t}(\alpha)$ for every $\alpha\in B_s$. Indeed, let $\{F_i \mid i < \omega\}$ be an infinite set of thread coverings. Using the finite $\Delta$-system lemma (the sunflower lemma), we may assume that $F_i \cap F_{i'} = F_\star$ for all $i \neq i'$. For every sufficiently large $\alpha$ which is not covered by $F_\star$, for all $i \neq i'$ and $f \in F_i \setminus F_\star, g \in F_{i'} \setminus F_\star$ $f\restriction \alpha \neq g \restriction \alpha$. For such $\alpha$, for every $i$ there is some $f \in F_i \setminus F_\star$ such that $f \restriction \alpha + 1 = e_\zeta$, $\zeta \in \name{t}(\alpha)$. Since all of them are distinct, the size of $\name{t}(\alpha)$ must be infinite, which is a contradiction.

Now if $s_0$ and $s_1$ are two Magidor-stems and $s_1$ is stronger than $s_0$ then every thread covering for $B_{s_1}$ is also a thread covering for $B_{s_0}$ since $B_{s_0}\subseteq B_{s_1}$.
In particular, there is a Magidor-stem $s_\star$ such that if $s$ is a stronger stem than $s_\star$ then the set of thread coverings for $B_s$ is identical with the set of thread coverings for $B_{s_\star}$.
Since $D'$ is expressible as a union of $B_s$s where every two Magidor-stems $s,s'$ are compatible, we conclude that in the generic extension there is a thread covering $F$ for $D'$.
It follows from Lemma \ref{lemcm} that $F\subseteq  a$, contradicting our assumption.
Hence the last case is covered, and we are done.

\hfill \qedref{thmmagforcing}

In order to accomplish the proof we need two lemmata.

\begin{lemma}
\label{lemcompleteness} The completeness lemma. \newline
The filter $\mathscr{G}$ is $\aleph_1$-complete in $V[G\ast H]$.
\end{lemma}

\par\noindent\emph{Proof}. \newline
Suppose that $(\zeta_n:n\in\omega)$ is a sequence of ordinals of $\mu^+$ in $V[G\ast H]$.
Let $A=\bigcap\{A_{\zeta_n}:n\in\omega\}$.
We need showing that $|A|=\mu^+$.
Fix an ordinal $\alpha_0\in\mu^+$.
We shall find an ordinal $\gamma\in\mu^+$ so that $\alpha_0<\gamma$ and $\gamma\in A$ in $V[G\ast H]$.

To begin with, choose a set $x\in V$ so that $|x|<\mu$ and it is forced by $\mathbb{P}\ast\mathbb{M}$ that $\{\zeta_n:n\in\omega\}\subseteq x$.
Such a set $x$ exists by virtue of Lemma \ref{lemcovering}.
Now choose an ordinal $\alpha_1\in\mu^+$ such that $\alpha_0<\alpha_1$ and moreover $\alpha_1>\ell g(e_{\zeta_n})$ for every $n\in\omega$.
By the genericity of the Cohen functions we can find an ordinal $\gamma\in\mu^+$ such that $\alpha_1\leq\gamma$ and $g_\alpha(\gamma)=1$ for every $\alpha\in x$.
It follows that $\gamma\in A_{\zeta_n}$ for every $n\in\omega$, namely $\gamma\in A$ in $V[G\ast H]$, so we are done.

\hfill \qedref{lemcompleteness}

Our second lemma provides us with the following property:

\begin{lemma}
\label{lemrqbeta} The evasion lemma. \newline
There are $(q,p)\in\mathbb{P}\ast\mathbb{M}$ and $\beta\in\mu^+$ such that:
\begin{enumerate}
\item [$(a)$] $\beta=M_i\cap\mu^+$ for some $i\in\mu^+$, and $\cf(\beta)=\mu$.
\item [$(b)$] $(q,p)\Vdash\check{\beta}\in\name{D}$.
\item [$(c)$] $(q,p)\Vdash\name{t}(\alpha)\cap R_{q\beta}=\varnothing$ for a $\mathscr{G}$-positive set of $\alpha\in D$.
\end{enumerate}
\end{lemma}

\par\noindent\emph{Proof}. \newline
We commence with a simple observation concerning the sets $R_{q\beta}$.
Suppose that $q_0\in\mathbb{P},\beta_0\in\mu^+$ and $A_0=\bigcap\{A_\xi:\xi\in R_{q_0\beta_0}\}$.
Since $|R_{q_0\beta_0}|<\mu$ we see that $A_0\in\mathscr{G}$.
Suppose that $q_0\leq q_1$ and $\beta_1\in A_0$.
We claim that $R_{q_0\beta_0}\cap R_{q_1\beta_1}=\varnothing$.
Indeed, if $\gamma\in R_{q_0\beta_0}\cap R_{q_1\beta_1}$ then $\beta_1\in A_\gamma$ since $\gamma\in R_{q_0\beta_0}$.
Hence $q_0\Vdash g_\gamma(\beta_1)=1$, so $q_1\Vdash g_\gamma(\beta_1)=1$ as well.
On the other hand, $\gamma\in R_{q_1\beta_1}$ so $q_1\Vdash g_\gamma(\beta_1)=0$, a contradiction.

Working in $V[G]$ we assume toward contradiction that the lemma fails.
Let $S=S^{\mu^+}_\mu$ and let $C=\{M_i\cap\mu^+:i\in\mu^+\}$, so $C\cap S$ is a stationary subset of $\mu^+$.
Every element of $C\cap S$ satisfies $(a)$, so for each $\beta\in C\cap S$ either $(b)$ fails or $(c)$.

Choose $(q_0,p_0)$ and $\beta_0\in C\cap S$ such that $(q_0,p_0)\Vdash\check{\beta}\in\name{D}$.
By our assumption toward contradiction we see that $(q,p)\Vdash\name{t}(\alpha)\cap R_{q_0\beta_0}\neq\varnothing$ for $\mathscr{G}$-almost every $\alpha\in D$.
Let $A_0=\bigcap\{A_\xi:\xi\in R_{q_0\beta_0}\}$.
Now by induction on $0<m\in\omega$ we choose $(q_m,p_m)$ and $\beta_m$ so that $(q_{m-1},p_{m-1})\leq(q_m,p_m), \beta_m\in C\cap S$ and $(q_m,p_m)\Vdash\beta_m\in\bigcap_{j<m}A_j\cap\name{D}$.

Pick up some $\alpha\in\bigcap_{m\in\omega}A_m\cap D$ such that $(q_m,p_m)\Vdash \name{t}(\alpha)\cap R_{q_m\beta_m}\neq\varnothing$ for every $m\in\omega$, and let $n=|\name{t}(\alpha)|$.
By the construction, $(q_{n+1},p_{n+1})\Vdash\name{t}(\alpha)\cap R_{q_j\beta_j}\neq\varnothing$ for every $j\leq n$.
By the observation at the beginning of the proof, if $j<\ell\leq n$ then $R_{q_j\beta_j}\cap R_{q_\ell\beta_\ell}=\varnothing$, so $(q_{n+1},p_{n+1})\Vdash|\name{t}(\alpha)|\geq n+1$, a contradiction.

\hfill \qedref{lemrqbeta}

The main result of this section points to the possibility that singular cardinals with countable cofinality behave in a different way than their friends with uncountable cofinality.
As mentioned in the introduction, it is easy to force a Dowker filter at $\mu^+$ even if $\mu>\cf(\mu)=\omega$ if we drop the requirement of strong limitude.

\begin{question}
\label{qcofomega} Is it consistent that $\mu$ is a strong limit singular cardinal with countable cofinality and there exists a Dowker filter over $\mu^+$?
\end{question}

\newpage 

\bibliographystyle{amsplain}
\bibliography{arlist}

\end{document}